\documentclass[a4paper, 12pt]{article}

\usepackage{graphicx}
\usepackage{amssymb}

\def\blankbox{{ ~\hfill$\rlap{$\sqcap$}\sqcup$}}

\begin{document}


\centerline {\bf Families of parameters for SRNT graphs}

\bigskip \bigskip

\centerline{Norman Biggs}

\bigskip \bigskip

\centerline{Department of Mathematics}

\centerline{London School of Economics}

\centerline{Houghton Street}

\centerline{London WC2A 2AE}

\centerline{U.K.}

\centerline{n.l.biggs@lse.ac.uk}

\bigskip

\centerline{October 2009}

\bigskip \bigskip

\centerline{\bf Abstract}
\bigskip

 The feasiblity conditions obtained in a previous report are refined, and used to
determine several infinite families of feasible parameters for strongly regular graphs with
no triangles.  The methods are also
used to improve the lower bound for the number of vertices, and to derive yet another
interpretation of the Krein bound.

\vfill \eject

 {\bf 1. Introduction}
\medskip

This paper is a continuation of my report on {\it Strongly Regular Graphs with No Triangles}, which will
be referred to as [SRNT1]. The graphs $X$ considered in that paper
are characterized by two parameters $k$ and $c$, according to the rules
\smallskip

{\leftskip 30pt

$\bullet$  $X$ is regular with degree $k$;

$\bullet$  any two adjacent vertices have no common neighbours;

$\bullet$  any two non-adjacent vertices have $c$ common neighbours.

\par}
\smallskip

It is convenient to rule out the pentagon and complete bipartite graphs,  so we impose the
conditions $k\ge 3$ and $k > c \ge 1$.
\smallskip

According to the standard theory,
the eigenvalues of the adjacency matrix of $X$ are
$k$ (with multiplicity 1) and the roots $\lambda_1, \lambda_2$
of the equation $\lambda^2 + c\lambda - (k-c) = 0$.  Furthermore, there is an integer $s >c$ such that
$c^2 + 4(k-c) =s^2$, where $s$ and $c$ have same parity, and the eigenvalues are the integers

$$ k = \frac {s^2 - c^2}{4} + c \qquad \lambda_1 = \frac{s-c}{2}, \qquad \lambda_2 = \frac{-s-c}{2}.$$

The multiplicities $m_1, m_2$ of $\lambda_1, \lambda_2$ are

$$m_1 = \frac{k}{2cs} \Big( (k-1+c)(s+c) - 2c\Big), \quad
  m_2 = \frac{k}{2cs} \Big( (k-1+c)(s-c) + 2c\Big).$$

\smallskip

In the appendix to [SRNT1] we took the basic parameters to be $q$ and $c$, where $q = \lambda_1$.
In that case $k$, $s$, and the number of vertices $n_q(c)$, are given by the formulae

$$k = c(q+1) +q^2, \quad  s = c+2q, \quad n_q(c) = Ac + B + D/c,$$

where $A = (q+1)(q+2)$, $B= 2q^3 + 3q^2 -q$, and $D = q^4 - q^2$.
\smallskip

Given $q$, it was established that $c$ must be in the range $1 \le c \le q(q+1)$, and must satisfy
\smallskip

{\leftskip 30pt
$c$ is a divisor of $q^4 - q^2$ \hskip 188pt (S1)
\smallskip

$c + 2q$ is a divisor of  $q^4 +3q^3 + 5q^2 +3q  +  q(q^4 - q^2)/c$ \qquad (S2).
\par}
\medskip

Using these results it was proved that  $n = n_q(c)$  is in  the range

$$ \Bigl\lceil   2q^3 + 3q^2 - q +2q(q+1) \sqrt{q^2 +q -2} \Bigr\rceil \;  \le  \; n \;  \le \;  q^2 (q+3)^2.$$

We shall now look more closely at the feasibility conditions, and use them to establish the existence of
infinite families of feasible parameters.  Also, in Section 5, we shall be able to
improve the lower bound for $n$.
\bigskip

\noindent{\bf 2.  Feasibility conditions revisited}
\medskip

It turns out that the feasibility conditions (S1) and (S2) can be stated more simply.
The formulae for $m_1$ and $m_2$ given above imply that

$$m_1 + m_2 =  \frac{k(k+c-1)}{c}, \qquad m_1 - m_2 = \frac{k(k+c-3)}{s}.$$

If these expressions are integers, then it follows $n= m_1 + m_2 + 1$ is an integer, but
 $m_1$ and $m_2$ may be half-integers. However, we can obtain a set of conditions equivalent to
(S1) and (S2) by adding the condition
$$m_1 + m_2 \equiv m_1 - m_2 \; ({\rm mod\;} 2).$$
\medskip

{\bf Lemma 1} \quad $m_1 + m_2$ is an integer if and only if
$$ c \quad {\rm  divides} \quad q^4 -q^2.$$
In that case, $m_1 + m_2$ and $(q^4 - q^2)/c$ have opposite parity.
\smallskip

{\it Proof} \quad  Using the fact that $k = c(q+1) +q^2$,
$$ \frac{k(k+c-1)}{c}  = c(q+1)(q+2) + (2q+1)(q^2 + q -1) + \frac{q^4 - q^2}{c}.$$
For all $q$, the first summand on the right-hand side is an even integer, and the second is an odd
integer.  Hence the result.
\blankbox
\medskip

{\bf Lemma 2 } \quad $m_1 - m_2$ is an integer if and only if
$$s = c+2q \quad {\rm  divides} \quad   q(q+1)(q+2)(q+3).$$
In that case $m_1 - m_2$ and $ q(q+1)(q+2)(q+3)/(c+2q)$ have opposite parity.
\smallskip

{\it Proof} \quad Putting $s = c+2q$,
$$ \frac{k(k+c-3)}{s}  = c(q+1)(q+2) - (3q^2 +7q +3) + \frac{q(q+1)(q+2)(q+3)}{c+2q}.$$
For all $q$, the first summand on the right-hand side is an even integer, and the second is an odd
integer.  Hence the result.
\blankbox

\vfill \eject

{\bf Theorem 1} \quad  The parameters $(q,c)$ are feasible for an SRNT graph if and only if $1 \le c \le q(q+1)$,
and the rational numbers

$$\alpha =  \frac{q^4 - q^2}{c} \qquad {\rm and}  \qquad \beta =  \frac{q(q+1)(q+2)(q+3)}{c+2q}$$

are integers having the same parity.

\blankbox

\bigskip

\noindent{\bf 3. Feasible parameters when $c$ is given}
\medskip

For a given $c$, the condition $c \le q(q+1)$ implies that

$$q \ge q_{min} = \Bigg\lceil \frac{\sqrt{4c+1} -1}{2}   \Bigg\rceil .$$

When $c$ is of the form $r(r+1)$, we have $q_{min} = r$, and these parameters are feasible for all $r$.
\smallskip

The next theorem shows that in almost all cases there is a corresponding number $q_{max}$.
\medskip

{\bf Theorem 2} \quad For all values of $c$ except $2,4,6$ there are only finitely many $q$ for which
$(q, c)$ is feasible.   Specifically, if $(q,c)$ is feasible then $q \le c(h-1)/2$,
where $h= |(c-2)(c-4)(c-6)|$, ($c \neq 2,4,6)$.
\smallskip

If $c$ is odd then  $q_m = c(h-1)/2$ is feasible.  If $c \ge 8$ is even, then
$q_m^* = c(h^*-1)/2$ is feasible, where $h^*$ is the integer $h/16$.
\smallskip

{\it Proof} \quad The following identity is easily verified:

$$c(c-2)(c-4)(c-6) - 16q(q+1)(q+2)(q+3) \qquad$$
$$ \qquad  = (c+2q)(c-2q-6)(c^2 - 6c+8 +12q +4q^2).$$
For a given value of $c$ it follows from the identity that if $\beta$  is an integer then
$q$ must be such that $c+2q$ divides $c(c-2)(c-4)(c-6)$. When $c$ is not one of $2,4,6$ this means that
$$ c+2q \; \le \; |c(c-2)(c-4)(c-6)|, \quad {\rm that \;is,} \quad q\le \frac{1}{2} c(h-1),$$
so the finiteness result is proved.
\smallskip

Suppose now that $c$ is odd and $q_m = c(h-1)/2$.  In this case $\alpha$ is even,
 and in order to show that $q_m$ is feasible, we must show that $\beta$ is also even.
We have $c+2q_m = ch$  and for $e=1,2,3$
$$ q_m +e  = (c-2e) X_e,  \qquad {\rm where} \qquad
      X_e = \frac{1}{2}\left( \frac{ch}{c-2e} -1 \right)  $$
is an integer.
It follows that
$$\beta = \frac{q_m(q_m+1)(q_m+2)(q_m +3)} {ch}  =  \frac{(h-1)X_1 X_2 X_3}{2}.$$
Since $h-1$ and at least one of $X_1,X_2, X_3$ is even, $\beta$ must also be even.
\smallskip

Now suppose that $c$  is even, $c =2c^*$, and $q_m^* = c(h^*-1)/2 = c^*(h^*-1)$.   We have
$ c+2q_m^* = 2c^*h^*$,  and for $e=1,2,3$
$$q_m^*  +e =  (c^*-e)X_e^* \qquad {\rm where} \qquad X_e^* = \frac{c^*h^*}{c^*-e} -1,$$
is an integer. It follows that
$$\beta = \frac{q_m^*(q_m^*+1)(q_m^*+2)(q_m^*+3)} {c^*h^*}  =  (h^*-1)X_1^* X_2^* X_3^*.$$
In this case $\alpha$ is even except when $c^* \equiv 2$ (mod $4$), and it can be checked that
$\beta$ has the same parity.  Hence $q_m^*$ is feasible when $c$ is even.
\blankbox
\medskip

For small even values of $c \ge 8$ it is easy to check by explicit computation that if $(q,c)$ is feasible
then $q \le q_m^*$, so that $q_m^*$ is the actual $q_{max}$.  It should be possible to
prove this in general.
\smallskip

The situation for $c=2,4,6$ is well
known, and is included here for completeness.
\medskip

{\bf Theorem 3} \quad The parameters $(q,c)$ are feasible for all $q$ when $c = 2,4,6$, except
in the cases $c=2$, $q \equiv 3$ (mod $4$), and $c=6$, $q \equiv 1$ (mod $4$).
\smallskip

{\it Proof} \quad Suppose $c=2$. Since $\alpha = (q^4 - q^2)/2$ is always an even integer, and $c+2q= 2(q+1)$,
we require that $\beta = q(q+2)(q+3)/2$ is also an even integer. It is easy to check that this holds except when $q \equiv 3$ (mod $4$).
\smallskip

The other cases are similar.
\blankbox
\medskip

On the basis of the preceding results, it is easy to list the possible values of $q$ for each $c$.
The values for small $c \neq 2,4,6$ are collected here for reference.
\smallskip

$c=1: \quad  q=1, 2, 7$.

$c=3: \quad  q=3$.

$c=5: \quad q=5$.

$c=7: \quad q=7,14,49$.

$c=8: \quad q=8$.

$c=9: \quad q=3,6,9,18, 27,48, 63, 90, 153, 468$.

$c=10: \quad q= 5, 10, 15, 25, 55$.

$c=11: \quad q= 11,12,22,33,44,77, 110,187,242,341,572,1727$.

$c=12: \quad q= 3, 4,6,9,12,14,24,30, 39, 54, 84,174$.

$c=13: \quad q= 13, 25, 39,52,65,130,208,403,494,637,1495,4498$.

$c=14: \quad q= 7,8,13,14,21,28,35,63,77,98,133,203,413$.
\smallskip

It will be noted that in many cases $q$ is a multiple of $c$. For example, $q=7c$ is feasible
for $c=1,2,4,6,7,9,11,12,14, \ldots \;$.  The following theorem covers several of these cases.
\medskip

{\bf Theorem 4} \quad  (1) \quad For any  positive integer $b$,  $q = bc$ is feasible whenever
$c \equiv 2,4,6$ modulo $2b+1$.
\smallskip

 (2) \quad For $b \equiv 1,7$ (mod $9$), $q = bc$ is feasible whenever $c \equiv 2,4,6$ modulo $(2b+1)/3$.
\smallskip

{\it Proof} \quad  (1) \quad  Clearly $ \alpha = b^2c(b^2c^2 - 1)$ is always an even integer.
Since $c + 2q = (2b+1)c$ we have

$$\beta = \frac {b(bc+1)(bc+2)(bc+3)}{2b+1}.$$

Let $c = (2b+1)f +2e$, where $e \in \{1,2,3\}.$  Then
$$bc +e = b((2b+1)f +2e) + e = (2b+1)(bf +e).$$
Thus $2b+1$ divides one of the factors in the numerator of $\beta$.   Since $2b+1$ is odd,
$\beta$ must be even.
\smallskip

(2) \quad As before, $\alpha$ is an even integer, and putting $t = (2b+1)/3$,
$$\beta = \frac {b(bc+1)(bc+2)(bc+3)}{3t}.$$
Let $c = tf +2e$, where $e \in \{1,2,3\}.$  Then
$$bc +e = b(tf +2e) + e = t(bf +3e).$$
Thus $t$ divides one of the factors in the numerator of $\beta$. Since $t$ is not divisible by $6$,
the quotient $(bc+1)(bc+2)(bc+3)/t$ is divisible by $6$, so $\beta$
must be an even integer.
\blankbox
\medskip

The first part of  theorem shows, for example, that $q=7c$ is feasible for
$c \equiv 2,4, 6$ (mod $15$), and the second part strengthens this to $c \equiv 2,4,6$ (mod $5$).
\smallskip

Theorem 4 can be reformulated in terms of the feasibility of $c = q/b$ for a given value of $q$.
This will be done in the next section.
\vfill \eject

{\bf 4. Feasible parameters when $q$ is given}
\medskip

We consider first  the factors of $q^4-q^2$ in the
ring ${\mathbb Z}[q]$. Since $q^4 -q^2 = (q-1) q^2 (q+1)$ and we require that
$c \le q(q+1)$, the relevant factors are
$$q-1, \; q, \; q+1, \; q^2 - q, \; q^2-1,  \; q^2, \; q^2 +q.$$
\smallskip

{\bf Theorem 5} \quad The parameters $(q,c)$ are feasible for all $q$ when
$c$ takes any one of the values $q, q^2 - q, q^2, q^2 + q$.
\smallskip

{\it Proof} \quad  When $c= q$, $\alpha = q^3 -q$ is always an even integer.
In this case $s = 3q$ and so we require that $\beta = (q+1)(q+2)(q+3)/3$ is even, which is clearly
true. The other cases are similar.
\blankbox

\bigskip

{\bf Theorem 6} \quad When $c= q-1, q+1, q^2 - 1$ the parameters $(q,c)$ are feasible for only a finite
set of values of $q$ in each case. Specifically,
\smallskip

$c= q-1$ is feasible only when $q= 2,5,7,12,47$.
\smallskip

$c= q+1$ is feasible only when $q= 1,3,13$.
\smallskip

$c= q^2 - 1$ is never feasible.
\smallskip

{\it Proof} \quad When $c =q-1$,  $\alpha = (q^4 - q^2)/c = q^3 + q^2$ is always an even integer.
In this case $s = 3q-1$, and we have the identity

$$  81q(q+1)(q+2)(q+3) = (3q-1)(27q^3 + 171q^2 + 354q + 280) + 280.$$
Hence, in order that $\beta = q(q+1)(q+2)(q+3)/(3q-1)$ should be an even integer, it is necessary that
$3q-1$ must evenly divide 280. The possibilities are $3q-1 = 5, 14, 20, 35, 140$, corresponding to
$q = 2, 5, 7, 12, 47$.
\smallskip

When $c=q+1$, $s= 3q+1$ and we have a similar identity

$$81q(q+1)(q+2)(q+3) = (3q+1)(27q^3 +153q^2 + 246q + 80) - 80.$$
Hence, in order that $\beta = q(q+1)(q+2)(q+3)/(3q+1)$ should be an even integer, it is necessary that
$3q+1$ must evenly divide 80. The possibilities are $3q+1 = 4, 10, 40$, corresponding to
$q = 1,3,13$.
\smallskip

When $c= q^2 -1$, $s= q^2 +2q -1$ and we have the identity

$$q(q+1)(q+2)(q+3) = (q^2 +2q -1)(q+2)^2 +  (2q+4).$$

So the condition is that  $q^2 + 2q -1$ evenly divides $2q+4$, which holds only
for the irrelevant values $q=1, c=0$.
\blankbox
\bigskip

Of course, there are other divisors of  $q^4 - q^2$, not covered by the arguments given above.
For example, if $c = q(q-1)/b$ is an integer, then $\alpha = bq(q+1)$ is an even integer.
Hence in this case $(q, c)$ is feasible provided that
$$\beta = \frac{b(q+1)(q+2)(q+3)}{q +2b-1},$$
is an even integer.  We have
$$b(q+1)(q+2)(q+3)   = (q +2b-1)Q + R,$$
where
$$Q= b(q^2 -(2b-7)q + +4b^2 -16b + 18), \qquad R= -4b(b-1)(b-2)(2b-3).$$
When $b=1$ and $b=2$, $Q$ is even and $R=0$, so the parameters are feasible for all $q$.
The case $c=q(q-1)$ is already covered in Theorem 5, but $c = q(q-1)/2$ is a new infinite family.
 For larger values of $b$
we get only finitely many feasible $q$ in each case.
\smallskip

We summarize the results so far by listing the values of $c$ that have been shown to be such that the parameters
$(q,c)$ are feasible for infinitely many $q$.
\medskip

{\leftskip 10pt
$c = 2$ for $q \equiv 0,1,2$ (mod $4$);
\smallskip

$c=4$ for all $q$;
\smallskip

$c=6$ for $q \equiv 0,2,3$ (mod $4$);
\smallskip

$c = q/b$ with $b\equiv 1,7$ (mod $9$) for $q \equiv 0$ (mod $b$) and $q \equiv 2b, 4b, 6b$ (mod $(2b+1)/3$);
\smallskip

$c=q/b$ with $b \equiv 0,2,3,4,5,6,8$ (mod $9$) for $q \equiv 0$ (mod $b$) and $q \equiv 2b,4b,6b$ (mod $2b+1$);
\smallskip

$c=q$ for all $q$;
\smallskip

$c= q(q-1)/2$ for all $q$;
\smallskip

$c=q(q-1)$ for all $q$;
\smallskip

$c=q^2$ for all $q$;
\smallskip

$c=q(q+1)$ for all $q$.
\par}

\vfill \eject

{\bf 5. The bounds for $n$}
\medskip

On the basis of the foregoing theory, it is possible to improve the lower bound for $n$ given in
the Appendix to [SRNT1].
\medskip

{\bf Theorem 7} \quad For all $q$ except $q=2,4,5,6,7,8, 12, 47$ the number
$n$ of vertices of an SRNT graph with $\lambda_1 = q$ lies in the range

$$4q^3 + 6q^2  \le n \le q^2(q+3)^2.$$

The bounds are attained by feasible parameters when $c = q$ and $c=q(q+1)$ respectively.
When $q= 2,5,7,12, 47$ the lower bound
is $4q^3 + 6q^2 - 2(q+1)$, and is attained when $c=q-1$.  When $q = 4,6,8$ the lower bound
is $4q^3 +6q^2 - 2(q-1) + 12/(q-2)$, and is attained when $c = q-2$.
\smallskip

{\it Proof} \quad  The lower bound given in [SRNT1] is obtained by showing that there is a unique
 minimum of $n_q(c) = Ac + B + D/c$, which occurs when
$$c = q \left( \frac{q-1}{q+2} \right) ^{\frac{1}{2}}.$$
Clearly this is just less than $q$, so we must examine values of $c$ in the neighbourhood of $q$.
We have shown in Section 4 that $c=q$ is feasible for all $q$, and in fact
$$n_q(q) = 4q^3 + 6q^2, \qquad n_q'(q) = 3q+1 > 0.$$
This implies that, if $c=q$ is not the actual minimum, then the minimum must occur at some feasible $c<q$.
\smallskip

Now, when $c=q-3$,
$$n_q(q-3) = 4q^3 + 6q^2 + 18 + 72/(q-3), \qquad n_q'(q-3) <0,$$
so values of $n_q(c)$ smaller than $4q^3 + 6q^2$ can only occur when $c=q-2$ or $c=q-1$.
\smallskip

When $c=q-2$
$$n_q(q-2) = 4q^3 + 6q^2 - 2q +2 + 12/(q-2),$$
which shows immediately that $q-2$ can only be feasible when  $q-2$ divides $12$, so we must check
the cases $q= 3,4,5,6,8,14$ individually. It turns out that only $q = 4,6,8$ are feasible.
\smallskip

The case $c= q-1$ was covered explicitly in Theorem 5, where we found that only the values
$2,5,7,12,47$ are feasible. In these cases we have  $n_q(q-1) = 4q^3 + 6q^2 - 2q -2$.
\blankbox \bigskip

{\bf 6.  Further properties of the second subconstituent}
\medskip

The results in [SRNT1] establish that, for an SRNT graph $X$, the second subconstituent
 $X_2$ is a connected graph of degree $k-c$ with diameter 2 or 3.
The only numbers that can be eigenvalues of $X_2$ are: $k-c$, $-c$ and the eigenvalues $\lambda_1$, $\lambda_2$
of $X$.   In terms of the parameters $(q,c)$, and in strictly decreasing order, these are
$$ q(q+c),  \quad q, \quad -c, \quad -(q+c).$$

Since $X_2$ is connected,  $q(q+c)$ has multiplicity $1$. Denote the multiplicities of $q,-c, -(q+c)$ by
$x,y,z$ respectively, and recall the standard formulae for $S_i =\sum \lambda^i$, $i= 0,1,2$.
In this case the resulting equations are
$$  S_0 = \ell, \qquad S_1 = 0, \qquad S_2 = \ell(k-c).$$
Putting $k= c(q+1) + q^2$, $\ell = k(k-1)/c$, these equations
 have a unique solution for $x,y,z$:

$$ x= \frac{(q+1)(q^2 + qc +c)(q^2 + 2qc + c^2 - 2c - q)}{c(c+2q)},$$

$$ y= q^2 + qc + c - 1, \hskip 142pt$$

$$ z= \frac{(q+c-1)(q^2 + qc +c)(q^2 + q - c)}{c(c+2q)}. \hskip 40pt$$

It can be verified that these values also satisfy the condition $S_3 = 0$, corresponding to the
fact that $X_2$ has no triangles.
\smallskip

Comparison with the formulae for the Krein parameters given in [SRNT1] shows that

$$qc(c+2q)z = (q^2 + qc + c)K_2.$$

Hence the fact that $z$ is a non-negative integer is equivalent to $K_2 \ge 0$. The fact that $z=0$
when $c = q(q+1)$ means that $X_2$ is an SRNT graph in this case.

\vfill \eject

{\bf References}
\medskip

[SRNT1] N.L. Biggs. Strongly regular graphs with no triangles. {\it Research Report} September 2009.
arXiv:0911.2160v1

\medskip

The references from [SRNT1] are repeated here for convenience.
\medskip

{\bf 1.} N.L. Biggs.  Automorphic graphs and the Krein condition. {\it Geom. Dedicata} (5) 1976 117-127.

{\bf 2.} A.E. Brouwer.  Strongly Regular Graphs. In: {\it Handbook of Combinatorial Designs}, ed. C. Colbourn,
J. Dinitz,  CRC Press, 1996.

{\bf 3.} A.E. Brouwer, A.M. Cohen, A.Neumaier.  {\it Distance-Regular Graphs}, Springer, Berlin 1989.

{\bf 4.} P.J. Cameron. Partial quadrangles. {\it Quart. J. Math. Oxford (2)} 26 (1975) 61-73.

{\bf 5} A.D. Gardiner. Antipodal covering graphs {\it J. Combinatorial Theory (Series B)} 16 (1974) 255-273.

{\bf 6.} C.D. Godsil. Problems in algebraic combinatorics. {\it Elect. J. Combinatorics} 2 (1995) F1.

{\bf 7.} C.D. Godsil, G. Royle. {\it Algebraic Graph Theory}, Springer, New York 2001.

{\bf 8.} D.M. Mesner. A new family of partially balanced incomplete block designs with some latin square
design properties. {\it Ann. Math. Statist.} 38 (1967) 571-581.

{\bf 9.} M.S. Shrikhande. Strongly regular graphs and quasi-symmetric designs. {\it Utilitas Mathematica}
3 (1973) 297-309.

{\bf 10.} S.S. Shrikhande. Strongly regular graphs containing strongly regular subgraphs.
{\it Proc. Indian Natl. Sci. Acad. Part A} 41 (1975) 195-203.

{\bf 11.} M.S.Smith. On rank 3 permutation groups.  {\it J. Algebra} 33 (1975) 22-42.

\end{document}